\documentclass[12pt]{article}
\usepackage{amsmath}
\usepackage{amssymb}
\newsymbol \blackbox 1004
\newcommand{\eh}{\hfill}\newlength{\sperr}

\newenvironment{proof}{{\settowidth{\sperr}{\bf\rm
Proof}%
\par\addvspace{0.3cm}\noindent\parbox[t]{1.3\sperr}
{\bf\rm P\eh r\eh o\eh o\eh f\eh }%
}}{\nopagebreak\mbox{}
$\blackbox$\par\addvspace{0.3cm}}

\def\ker{{\rm Ker\ }}
\newtheorem{Pa}{Paper}[section]
\newtheorem{Tm}[Pa]{{\bf Theorem}}

\newtheorem{Cy}[Pa]{{\bf Corollary}}
\newtheorem{Rk}[Pa]{{\bf Remark}}
\newtheorem{Pb}[Pa]{{\bf Problem}}
\newtheorem{Ee}[Pa]{{\bf Example}}

\newtheorem{Pn}[Pa]{{\bf Proposition}}
\newcommand{\CC}
{{\mathchoice {\setbox0=\hbox{$\displaystyle\rm
C$}\hbox{\hbox
to0pt{\kern0.4\wd0\vrule height0.9\ht0\hss}\box0}}
{\setbox0=\hbox{$\textstyle\rm C$}\hbox{\hbox
to0pt{\kern0.4\wd0\vrule height0.9\ht0\hss}\box0}}
{\setbox0=\hbox{$\scriptstyle\rm C$}\hbox{\hbox
to0pt{\kern0.4\wd0\vrule height0.9\ht0\hss}\box0}}
{\setbox0=\hbox{$\scriptscriptstyle\rm C$}\hbox{\hbox
to0pt{\kern0.4\wd0\vrule height0.9\ht0\hss}\box0}}}}

\title{Operator Bezoutiant and Roots of Entire Functions, Concrete Examples}

\author{L.A. Sakhnovich}

\date{}
\begin{document}
\maketitle

\begin{abstract}
In this paper we use the Bezoutiant method to describe  the conditions under which
two entire functions have not common roots. We apply the general results to concrete
examples. In particular we consider the Bessel functions.

\end{abstract}

{MSC(2000) Primary 30D20; Secondary 30D35, 47B35.}
\\

\\

735 Crawford ave., Brooklyn, 11223, New York, USA.\\
 E-mail address: lev.sakhnovich@verizon.net

\maketitle

 \section{Introduction}\label{Intro}
\setcounter{equation}{0}
 The matrix Bezoutiant  is used in order to define the number of common zeroes of two polynomials $f(z)$ and $g(z)$. M.G. Krein extended the notion of Bezoutiant to entire functions of the form
 \begin{equation}\label{1.1}
 F(z)=1+\int_{0}^{\omega}e^{izt}\overline{\Phi(t)}dt,\quad \Phi(t){\in}L(0,\omega).\end{equation}
 The result by M.G.Krein was not published and I became acquainted with it from the manuscript given to me by M.G.Krein in 1974. In 1976  I. Gohberg and G. Heinig published the article [4], in which deduced Krein's theorem and generalized it for the matrix functions
$F(z)$ of type (1.1). In the same 1976 we extended the Krein's theorem to the class of functions of the form [7]:
\begin{equation}
F(z)=1+iz\int_{0}^{\omega}e^{izt}\overline{\Phi(t)}dt,\quad \Phi(t){\in}L(0,\omega).\end{equation} Later  in the Bezoutiant theory a number of important and interesting results was published (see[3],[5]). In particular these results established
the connection between the two following problems:
\begin{Pb}\label{Pbm}  To find  the number of common zeroes of  the two entire functions $F_{1}(z)$ and $F_{2}(z)$.\end{Pb}
\begin{Pb}\label{Pbd}To describe the dimension of the Bezoutiant kernel.\end{Pb}
 The  Problem 1.2 is solved with the help of the finite number of  arithmetic actions which provides the effectiveness of the Bezoutiant approach when $F_{1}(z)$ and $F_{2}(z)$ are  polynomials and  the corresponding Bezoutiant is a matrix. In the operator case the situation is more complex. Up till now there hasn't
been a single concrete example of effective application of the operator Bezoutiant theory.
The main aim of this work is the construction of such examples.
We apply the operator Bezoutiant theory to the entire functions of the form
\begin{equation}
F_{k}(z)=\int_{0}^{a}e^{izt}\overline{\Psi_{k}(t)}dt.\end{equation}
We investigate in detail a class $Z$ of the functions $F_{k}(z)$ of form (1.3)
when $\Psi_{k}(t)$ is a polynomial with algebraic coefficients.
We  proved the following assertion:
\begin{Tm}\label {TmIbp} Let the following conditions be fulfilled.\\
1.The functions  $F_{k}(z)$ have the form (1.3).\\
2. \begin{equation}\Psi_{1}(x){\ne}\overline{\Psi_{2}(a-x)}\end{equation}
3.\begin{equation}\int_{0}^{a}\Psi_{k}(x)dx{\ne}0,\quad k=1,2.\end{equation}
Then  the corresponding functions $F_{1}(z)$  and $\overline{F_{2}(\overline{z})}$ haven't common zeroes.If $\Psi_{1}(x){\ne}\overline{\Psi_{1}(a-x)}$, then the corresponding function $F_{1}(z)$ hasn't real zeroes and hasn't conjugate pairs of zeroes. \end{Tm}
We shall use the following equality
\begin{equation}\overline{F_{2}(\overline{z})}=\int_{0}^{a}e^{-izt}\Psi_{2}(t)dt=
e^{-iaz}\int_{0}^{a}e^{izt}\Psi_{2}(a-t)dt.\end{equation}
Hence the next assertion is true.
\begin{Pn}\label{PnF}The functions $\overline{F_{2}(\overline{z})}$ and
\begin{equation}F_{2,1}(z)=\int_{0}^{a}e^{izt}\Psi_{2}(a-t)dt\end{equation}
have the same zeroes.\end{Pn}
\begin{Rk}It is important , that the function $F_{2,1}(z)$ belongs to the class Z.\end{Rk}
\begin{Ee}\label{Eet}Let $\Psi(t)=t^{n}$, where $n{\geq}0$ and integer. In this case we have (see [2]):
\begin{equation}
F(n,z)=-(-i)^{n+1}\frac{d^{n}}{dx^{n}}[x^{-1}(1-\mathrm{cos}x-i\mathrm{sin}x)]{\in}Z.
 \end{equation}\end{Ee}
 \begin{Cy}\label{Cyff}The different functions $F(n_{1},z)$ and $F(n_{2},z))$ ,defined by (1.8)  haven't common zeroes.\end{Cy}
\begin{Ee}\label{Eeta}
Let $\Psi(t)=t^{n}(a-t)^{m}$, where $n$  and $m$ are integer and $n{\geq}0$, $m{\geq}0$.
 The corresponding function $F(n,m,z)$ belongs to the class $Z$. If $n=m$ we have
 \begin{equation}F(n,n,z)=\sqrt{\pi}\Gamma(z/2)^{-(n+1/2)}J_{(n+1/2)}(z){\in}Z,\end{equation}
 where $\Gamma(z)$ is Euler Gamma function, $J_{\nu}(z)$ is Bessel function.
The functions $J_{\nu}(z)$ form a subclass $Z_{1}$ of class $Z$.\\
 For subclass $Z_{1}$
the Theorem 1.3  has been well-known more than a hundred years [6],[14].
\end{Ee}
\begin{Ee}\label{Eeop}\emph{Open problem.}\\It is interesting to use our approach to the case \begin{equation}\Psi(t)=t^{n+1/2}(a-t)^{m+1/2},\end{equation} where $n$  and $m$ are integer and $n{\geq}0$, $m{\geq}0$. The results of sections 2-5 are true for case (1.10) too.\end{Ee}
The class of Bessel functions $J_{(n)}(z),$ where $n$  is  integer and $n{\geq}0$  can be reduced to this case (see (1.9)). Let us formulate the famous Bourget's hypothesis [14]:\\
\textbf{Bourget's hypothesis}\emph{Two functions $J_{n}(z)$ and $J_{m}(z)$ ,where $n{\ne}m$, haven't common zeroes other than the origin.}\\ With the help of the  Siegel theorem [12],[14] it was proved that the
Bourget's hypothesis is true.\\Now we formulate the generalized Bourget's hypothesis\\
\textbf{Generalized Hypothesis }\emph{Let one of the following inequalities is true:\\
Either $n_{1}{\ne}m_{2}$ or $n_{2}{\ne}m_{1}.$
Two  functions $F_{1}(z)$ and  $F_{2,1}(z)$  haven't common zeroes.\\ If $n_{1}{\ne}m_{1}$,then the corresponding function $F_{1}(z)$ hasn't real zeroes and hasn't conjugate pairs of zeroes.}\begin{Rk}The functions $F(z)$ from Examples 1.8
and 1.9 can be expressed in the terms of the confluent hypergeometric function $\Phi(b,c,z)$, if we use the representation [1]:
\begin{equation}\Phi(b,c,z)=\frac{\Gamma(c)}{\Gamma(b)\Gamma(c-b)}\int_{0}^{1}e^{zt}t^{b-1}(1-t)^{c-b-1}dt,\end{equation}
where $Re{c}>Re{b}>0.$\end{Rk}
Let us briefly explain the structure of the paper.\\
In sections 1-4 we describe the main notions and results from paper [7] (see [8], Ch.5).
In section 5 we construct  for  case (2.3)  the  operator  Bezoutiant $T$ in the explicit form. By this  construction we use the methods of the operator identities
[5],[11]. This result gives an effective analytic method for solving the formulated problems 1.1 and 1.2. Section 6 is contained  the proof of Theorem 1.3, and the investigation of    Examples 1.6 and 1.8.
\section{Main notions }\label{Main}
\setcounter{equation}{0}
By $[H_{1},H_{2}]$ we denote the set of linear bounded operators acting from the Hilbert space $H_{1}$  into the Hilbert space $H_{2}$. The space of constant $m{\times}1$ vectors we denote by G .Now we introduce the $m{\times}m$ matrix
functions
\begin{equation}F_{1}(z)=I_{m}-zP^{\star}(I-Az)^{-1}\Pi,\end{equation}
\begin{equation}  F_{2}(z)=I_{m}-zQ^{\star}(I-Az)^{-1}\Pi.\end{equation}
Here the operators A, P , Q and $\Pi$ are such  that
$A{\in}[L_{m}^{2}(0,\omega),L_{m}^{2}(0,\omega)]$, $\Pi{\in}[G,L_{m}^{2}(0,\omega)]$,
$P^{\star}{\in}[L_{m}^{2}(0,\omega),G]$, $Q^{\star}{\in}[L_{m}^{2}(0,\omega),G]$.\\
 Let us note that the representation of the given matrix functions $F_{1}(z)$ and $F_{2}(z)$ is called the realization. The methods of realization
are well-known\\ ( see [9]).\\
 \emph{Further we assume that the spectrum of the operator A  coincides with zero.}\\
 Hence the functions $F_{1}(z)$ and $F_{2}(z)$ defined respectively by (2.1) and (2.2) are  entire matrix  functions.
 Let us associate with pair $F_{1}(z)$ and $F_{2}(z)$ the operator identity
 \begin{equation}TB-C^{\star}T=N_{2}N_{1}^{\star}, \end{equation} where
  \begin{equation}B=A+{\Pi}P^{\star},\quad C=A+{\Pi}Q^{\star}, \quad N_{1}^{\star}={\Pi^{\star}}T.\end{equation}
The operators $B,C,T$ and $N_{1},N_{2}$ are such that $B,C{\in}[L_{m}^{2}(0,a),L_{m}^{2}(0,a)]$, $T{\in}[L_{m}^{2}(0,a),L_{m}^{2}(0,a)]$ and $N_{1},N_{2}{\in}[G,L_{m}^{2}(0,a)].$  We want to stress that\\
 $N_{1}^{\star},N_{2}^{\star}{\in}[L_{m}^{2}(0,a),G];$  $\Pi^{\star}{\in}[L_{m}^{2}(0,a),G]$.\\
By $L_{T}$ we denote the kernel of $T$ , by $L_{1}$ we denote the maximal invariant subspace in respect to $B$ such that
\begin{equation}\label{1.5 }N_{1}^{\star}L_{1}=0.\end{equation}
In paper [7] (see [8], Ch.5) we proved the following assertion.
\begin{Tm}\label{Tmdim}Let the following conditions be fulfilled :\\
1)Relations (2.3), (2.4) are true.\\
2)If $M$ is an invariant subspace in respect to $A^{\star}$ and if $\Phi^{\star}M=0$ then
$M=0$.\\
In this case the equality
\begin{equation} L_{1}=L_{T} \end{equation}
is valid.\end{Tm}

\begin{proof}.
Equalities   (2.3), (2.4) imply that
\begin{equation}L_{T}{\in}\ker{N_{1}^{\star}}.\end{equation}
It follows from relation (2.3) that
\begin{equation}TBf=0,\quad f{\in}L_{T}, \end{equation}i.e. the subspace $L_{T}$ is $B$ invariant. Hence in view of (2.3) we have
\begin{equation}
L_{T}{\in}L_{1} .\end{equation} Operator identity (2.3) implies that
the subspace $H_{1}=\overline{TL_{1}}$ is $C^{\star}$ invariant. Due to (2.5) the relation $\Pi^{\star}H_{1}=0$ is valid. It means that on the subspace $H_{1}$ the operators $C^{\star}$ and $A^{\star}$ coincide. Using condition 2) of the theorem we deduce the equality $H_{1}=0$, i.e.
\begin{equation}
L_{1}{\in}L_{T} .\end{equation}The assertion of the theorem follows directly from (2.9) and (2.10).
\end{proof}

\begin{Ee}\label{EeO} Let us consider the extreme case, when  $T=0$. In view of relation (2.4) we have $N_{1}=0$.It means that $L_{1}=L_{T}=L_{m}^{2}(0,\omega)$.\end{Ee}

\begin{Ee}\label{EeA} Let us consider another extreme case, when $P=Q=0$ and
\begin{equation}Af=i\int_{0}^{x}f(t)dt,\quad f(x){\in}L^{2}(0,a).\end{equation}\end{Ee} In this case we have
\begin{equation}(A-A^{\star})f=i\int_{0}^{a}f(t)dt.\end{equation}It follows from (2.12) that
\begin{equation}T=I,\quad N_{1}g=g,\quad N_{2}g=ig,\quad g{\in}G_{1}.\end{equation}We see that $L_{T}=0$.It is well-known ([13],Ch.11) that the operator $A$,defined by relation (2.11), hasn't invariant subspaces orthogonal to 1. Hence $L_{1}=0$, i.e. we have again the equality $L_{1}=L_{T}$.
\section{Properties of the operator B}\label{Prop}
\setcounter{equation}{0}
\emph{Further we consider only the case when $\mathrm{dim}{G}=1$.}\\
In this section we formulate the well-known properties of the operator $B$ (see [7],[5]).
\begin{Pn}If $z$ is a regular point of $(I-Az)^{-1}$ and $F^{-1}_{1}(z)$ then $z$ is a regular point of  $(I-Bz)^{-1}$ and
\begin{equation} (I-Bz)^{-1}=(I-Az)^{-1}+z(I-Az)^{-1}{\Pi}F^{-1}_{1}(z)P^{\star}(I-Az)^{-1}
\end{equation}\end{Pn}
\begin{Rk}\label{Rkz}In view of (2.4) and (3.1) we have
\begin{equation} (I-Bz)^{-1}\Pi=(I-Az)^{-1}{\Pi}F^{-1}_{1}(z).\end{equation}\end{Rk}
\begin{Pn}The following relation \begin{equation}
 (B-zI)^{p+1}=\sum_{s=0}^{p}(A-zI)^{p-s}{\Pi}P^{\star}(B-zI)^{s}+ (A-zI)^{p+1}\end{equation}
is true.\end{Pn}
 Let $\lambda$ be an eigenvalue of operator $B$ and let $f_{p}$ be a corresponding
 root vector , i.e.
 \begin{equation} (B-{\lambda}I)^{p+1}f_{p}=0,\quad (B-{\lambda}I)^{p}f_{p}{\ne}0.\end{equation} Equation (3.3) implies that
\begin{equation}f_{p}=\sum_{s=0}^{p}(A-{\lambda}I)^{-s-1}h_{s},\end{equation}
where
\begin{equation}h_{s}=-{\Pi}P^{\star}(B-{\lambda})^{s}f_{p}.
\end{equation}Let us now consider the chain of the root vectors
\begin{equation}f_{p-k}=(B-{\lambda})^{k}f_{p},\quad 0<k{\leq}p.\end{equation}
It follows from (3.3) and (3.7) that
\begin{equation}f_{p-k}=\sum_{s=0}^{p-k}(A-{\lambda}I)^{-s-1}h_{s+k}.\end{equation}
In view of (3.8) we have
\begin{equation}f_{0}=(A-{\lambda}I)^{-1}h_{p}\end{equation}
where
\begin{equation}h_{p}=-{\Pi}P^{\star}f_{0}.\end{equation}
\begin{equation}(I-Bz)^{-1}\Pi=(I-Az)^{-1}{\Pi}F^{-1}_{1}(z).\end{equation}
\begin{Pn}\label{Pne}If the operators $A$ and $B$ do not have common eigenvalues then
\begin{equation}P^{\star}f_{0}{\ne}0.\end{equation}\end{Pn}
Let $\mu$ be an eigenvalue of the operator $C$ and let $g_{q}$ be the root vector of the order q. The following statement is true .
\begin{Pn}\label{Pnq}If the operators $A$ and $C$ do not have common eigenvalues then
\begin{equation} Q^{\star}g_{0}{\ne}0.\end{equation}\end{Pn}
\section{The explicit form of  Bezoutiant}\label{Expl}
\setcounter{equation}{0}
In this section we construct the operator Bezoutiant $T$ in the explicit form.
Let us consider the entire functions
\begin{equation}F_{k}(z)=\int_{0}^{a}e^{izt}\overline{\Psi_{k}(t)}dt,\quad,(k=1,2),\quad \Psi_{k}(t){\in}L(0,a).\end{equation}
From relation (4.1) we obtain that
\begin{equation}F_{k}(z)=1+iz\int_{0}^{a}e^{izt}\overline{\Phi_{k}(t)}dt,\end{equation}
where
\begin{equation}\Phi_{k}(t)=\frac{1}{R_{k}}\int_{t}^{a}\Psi_{k}(s)ds,\quad R_{k}=\int_{0}^{a}\Psi_{k}(u)du.\end{equation} Further we suppose that
\begin{equation}R_{k}{\ne}0,\quad k=1,2.\end{equation}
Formula (4.2) can be represented in the form
\begin{equation}F_{k}(z)=1-zP_{k}^{\star}(I-Az)^{-1}1,\end{equation}
where the operator $A$ is defined by relation (2.12) and
\begin{equation}P_{k}^{\star}f=-i\int_{0}^{a}f(t)\overline{\Phi_{k}(t)}dt.\end{equation}We use here the equality
\begin{equation}(I-Az)^{-1}1=e^{izx}.\end{equation}
We choose $\alpha$ and $\beta$ so that $\overline{\alpha}+\beta{\ne}0$ and put
\begin{equation}M_{1}(x)=\Phi_{2}(x)-\beta M_{2}(x),\quad M_{2}(x)=[\Phi_{2}(x)+\overline{\Phi_{1}(a-x)}-1]/(\overline{\alpha}+\beta).\end{equation}
To the pair of functions $F_{1}(z)$ and  $F_{2}(z)$ we assign the operator $T$
acting in $L^{2}(0,a)$ and defined by formulas (see [7]):
\begin{equation}Tf=\frac{d}{dx}\int_{0}^{a}f(t)\frac{\partial}{{\partial}t}\Phi(x,t)dt,
\end{equation} where
\begin{equation}\Phi(x,t)=\frac{1}{2}\int_{x+t}^{2a-|x-t|}Q(\frac{s+x-t}{2},\frac{s-x+t}{2})ds,
\end{equation}
\begin{equation}Q(x,t)=M_{2}(a-t)M_{1}(x)+[1-M_{1}(a-t)]M_{2}(x).\end{equation}
We introduce the matrices
  \begin{equation}A(x)=[M_{2}(x),1-M_{1}(x)],\quad B(x)=\mathrm{col}[M_{1}(x),M_{2}(x)].\end{equation} It follows from (4.3),(4.8) and (4.12) that
  $A(0)=0,\quad B(2)=0.$ Using formulas (4.9)-(4.12) we represent the operator $T$ in the form
  \begin{equation}Tf=c\int_{0}^{a}f(t)U(x,t)dt,\quad c=-\frac{1}{R_{1}R_{2}(\overline{\alpha}+\beta)}{\ne}0,\end{equation}
  where
  \begin{equation}U(x,t)=\int_{t}^{a}[\Psi_{2}(a-s)\overline{\Psi_{1}(a-s-x+t)}-
  \Psi_{2}(s+x-t)\overline{\Psi_{1}(s)}]ds,\end{equation}when  $(x<t)$ and
   \begin{equation}U(x,t)=\int_{t}^{a+t-x}[\Psi_{2}(a-s)\overline{\Psi_{1}(a-s-x+t)}-
  \Psi_{2}(s+x-t)\overline{\Psi_{1}(s)}]ds,\end{equation}when  $(x>t)$.
 \begin{Pn}\label{Bound} Let the condition $\Psi_{k}(x){\in}L(0,a)$ (k=1,2) be fulfilled.Then the operator $T$ defined by formulas $(4.13) - (4.15)$ is bounded in the space $L^{2}(0,a)$\end{Pn}

\begin{proof}.
Using formula
\begin{equation} \Psi_{k}(s)=0,\quad s{\notin}[0,a] \end{equation}
we extend the functions $\Psi_{k}(s)$.
It follows from (4.14) and (4.15)
   that \\$|U(x,t)|{\leq}h(x-t)$ , where
   \begin{equation}h(x)=\int_{0}^{a}[|\Psi_{2}(a-s)\overline{\Psi_{1}(a-s-x)}|+
  |\Psi_{2}(s+x)\overline{\Psi_{1}(s)}|]ds,\quad |x|{\leq}a.\end{equation}It is  easy to see that
   \begin{equation}\int_{-a}^{a}h(x)dx<\infty.\end{equation}
Hence the operator $T$ is bounded.
   The proposition is proved.
   \end{proof}
   
   In paper [7] the  following relations are
deduced:
\begin{equation}TB_{1}-B_{2}^{\star}T=N_{2}N_{1}^{\star}, \end{equation} where
  \begin{equation}B_{k}=A+{\Pi}P_{k}^{\star},\quad N_{2}g=-i(\overline{\alpha}+\beta)M_{2}(x)g, \quad N_{1}g=\overline{M_{2}(a-x)}g.\end{equation}
  A direct calculation shows that
  \begin{equation}T^{\star}1=\overline{M_{2}(a-x)}.\end{equation}  Relation  (4.21) can be written in the form (see (2.4)):
  \begin{equation}
  N_{1}^{\star}=\Pi^{\star}T.\end{equation}
   Now let us consider the function
  \begin{equation} F_{2,1}(z)=\overline{F_{2}(\overline{z})}e^{iaz}=1+iz\int_{0}^{a}e^{izt}\overline{\Phi_{2,1}(t)}dt,\end{equation}
where
\begin{equation} \Phi_{2,1}(t)=1-\overline{\Phi_{2}(a-t)}.\end{equation}
 It follows  from relations (4.8) and  (4.24) that
  \begin{equation}\Phi_{1}(t)-\Phi_{2,1}(t)=(\alpha+\overline{\beta})\overline{M_{2}(a-x)}
  \end{equation}  Taking into account relations (4.19) and (4.20)we deduce that the functions $F_{1}(z)$ and $F_{2,1}(z)$ satisfy all the conditions of Theorem 4.2.
  We note that the zeros of functions $\overline{F_{2}(\overline{z})}$ and $F_{2,1}(z)$ coincide. Hence the following statement is true.
\begin{Tm}\label{TmN}Let the condition $\Psi_{k}(x){\in}L(0,a)$ be fulfilled and $\mathrm{dim}L_{T}=N<\infty,$  where the operator $T$ is  defined by formulas (4.13) - (4.15).Then the number of common zeroes of $F_{1}(z)$ and $F_{2,1}(z)$ is equal to N as well.\end{Tm}
\begin{Rk}\label{RkT}It is important that the operator $T$ is constructed in the terms
of the given functions $F_{1}(z)$ and $F_{2}(z)$,i.e in terms of $\Psi_{1}(x)$
and $\Psi_{2}(x)$.\end{Rk}  \section{Classes of entire functions without common zeroes}\label{Classes}
  \setcounter{equation}{0}
\begin{Ee}\label{Eeone}Let the functions $\Psi_{k}^{(p)}(x) \quad (k=1,2;\quad 0{\leq}p{\leq}Q+1)$ be continuous.
 Then the relation
 \begin{equation} \frac{d^{Q+1}}{dx^{Q+1}}(Tf)=L(D)f(x)+\int_{0}^{a}f(t)[V(x-t)+W(x,t)]dt \end{equation} is true.\end{Ee}
 Here the kernel $W(x,t)$ is continuous and the kernel $V(x-t)$ and the differential operator $L(D)$ are defined by the relations
 \begin{equation} V(u)=\sum_{p+k=Q}[(-1)^{k+1}\Psi_{2}^{(p)}(u)\overline{\Psi_{1}^{(k)}(0)}+
 \Psi_{2}^{(k)}(a)\overline{\Psi_{1}^{(p)}(a-u)(-1)^{(p)}}],\end{equation}
 \begin{equation} L(D)=\sum_{p+k+s=Q-1}[(-1)^{k+1}\Psi_{2}^{(p)}(0)\overline{\Psi_{1}^{(k)}(0)}+
 \Psi_{2}^{(k)}(a)\overline{\Psi_{1}^{(p)}(a)(-1)^{(p)}}]D^{s}.\end{equation}
We denote by $D$ the operator $D=\frac{d}{dx}$ and  by $r$ the order of the
  differential operator $L(D)$ defined by relation(5.3). If $r=0$ then
  $L(D)f(x)={\alpha}f(x),\quad \alpha{\ne}0.$
In view of (5.1) the following assertion is true.
\begin{Pn}\label{Pnr} If the order $r$ of the differential operator $L(D)$ is non-negative, then $\mathrm{dim}L_{T}<\infty$ and the number of common zeroes of the corresponding functions  $F_{1}(z)$ and $F_{2,1}(z)$  is equal to $N =\mathrm{dim}L_{T}<\infty.$\end{Pn}
\begin{Ee}\label{EeT}Let us consider the case when
\begin{equation}\Psi_{1}(x)=\overline{\Psi_{2}(a-x)}\end{equation}\end{Ee}
In this case we have
\begin{equation}\overline{F_{2}(\overline{z})}=e^{-iza}F_{1}(z).\end{equation}
Using relations (5.13)-(5.15) and (5.5) we obtain the following assertion.
\begin{Pn}\label{PnTo}If relation (5.4) is true,then all the zeroes of the corresponding functions  $F_{1}(z)$ and $F_{2,1}(z)$ coincide and $T=0,\quad L_{T}=L^{2}(0,a).$\end{Pn}
\begin{Ee}\label{EePol} Let us consider the important special case when
 \begin{equation}\Psi_{k}(x)=\sum_{p=0}^{Q_{k}}b_{k,p}x^{p},\quad b_{k,Q_{k}}{\ne}0.
 \end{equation}\end{Ee}
 Now we shall formulate and prove the main theorem of this section.
 \begin{Tm}\label{TmP}Let the following conditions be fulfilled.\\
1.The functions  $\Psi_{k}(x)$ have the form (5.6),where
\begin{equation}Q=Q_{1}{\geq}Q_{2}.\end{equation}
2. \begin{equation}\Psi_{1}(x){\ne}\overline{\Psi_{2}(a-x)}\end{equation}
3.\begin{equation}\int_{0}^{a}\Psi_{k}(x)dx{\ne}0,\quad k=1,2.\end{equation} Without loss of generality we shall suppose that
\begin{equation}\int_{0}^{a}\Psi_{k}(x)dx=1,\quad k=1,2.\end{equation}
4.The numbers $a$ and $b_{p,k}$ are algebraic.\\
5.The order $r$ of the corresponding differential operator $L(D)$ is non-negative.
Then  the corresponding functions $F_{1}(z)$  and $F_{2,1}(z)$ haven't common zeroes.If $\Psi_{1}(x){\ne}\overline{\Psi_{1}(a-x)}$, then the corresponding function $F_{1}(z)$ hasn't real zeroes and hasn't conjugate pairs of zeroes. \end{Tm}

\begin{proof}.
 It  follows from paper [7] that there exists such $z_{j}$ that $Tf_{j}=0$, where $f_{j}=e^{z_{j}}$.  Hence we have (see(5.1)):
\begin{equation} L(z_{j})f_{j}(x)+\int_{0}^{x}f_{j}(t)V(x-t)dt=0 . \end{equation}
As the Volterra operator $T_{1}f=\int_{0}^{x}f(t)V( x-t)dt$ can not have the eigenvalues different from the zero therefore
\begin{equation} D (z_{j})=0.\end{equation}
Using  relations (5.11), (5.12) and Titchmarsh's theorem (see [13],Ch.11) we deduce that
\begin{equation} V(u){\equiv}0.\end{equation}
Now we write  the following equality (see[2]):
\begin{equation}F(m,z)=\int_{0}^{a}e^{itz}t^{m}dt=
-(-i)^{m+1}\frac{d^{m}}{dx^{m}}[x^{-1}(1-\mathrm{cos}ax-i\mathrm{sin}ax)].\end{equation}
The functions $F_{1}(z)$ can be represented in the form
 \begin{equation}F_{1}(z)=P(z)\mathrm{cos}az+Q(z)\mathrm{sin}(az)+R(z),\end{equation}
where $P(z),\quad Q(z)$ and $(z)$ are rational functions with algebraic coefficients. The equation $F_{1}(z)=0$ is equivalent to the equation
 \begin{equation}P(z)(1-t^{2})+2Q(z)t+R(z)(1+t^{2})=0,\end{equation}
where $t=\mathrm{tg}az/2$. According to relation (5.12) the common zero $z_{j}$ of the equations
$F_{1}(z)=0$ and $F_{2,1}(z)$ is an algebraic number. Relation (5.16) implies that $t=\mathrm{tg}az_{j}/2$ is algebraic number too. This fact contradicts to the following well-known assertion (see [12],[14]):\\
\emph{If $z_{j}$ is an algebraic number then  $\mathrm{tg}az_{j}/2$ is a transcendental number}.\\
 Hence the assertion of the theorem is true.
 \end{proof}
 
 \begin{Pn}\label{PnC}Let conditions 1-4 of Theorem 5.6 be fulfilled. Then $\mathrm{dim}L_{T}=0$
 and the corresponding functions $F_{1}(z)$  and $F_{2}(z)$ haven't common zeroes.If $\Psi_{1}(x){\ne}\overline{\Psi_{1}(a-x)}$, then the corresponding function $F_{1}(z)$ hasn't real zeroes and hasn't conjugate pairs of zeroes. \end{Pn}
 
 \begin{proof}.
 We assume that condition 5 of Theorem 5.6  is not valid, i.e $L(D)=0.$ If  $\mathrm{dim}L_{T}>0$,then according to Titchmarsh theorem (see [13],Ch.11) we have $V(u){\equiv}0$. It follows from (5.11) that\\
 $\frac{d^{Q+1}}{dx^{Q+1}}(Tf)=0,$i.e. $Tf$ is a polynomial in respect to $x$ of the order $Q$. Now we use the relations $T1=M_{2}(x)$, where  $M_{2}(x)$ is a polynomial
 of the order $P{\leq}Q.$ Then we obtain
$$\frac{d^{Q+1}}{dx^{Q+1}}(TB_{1}-B_{2}^{\star}T)=i\frac{d^{Q}}{dx^{Q}}T=0.$$
From the last relation we deduce that $M_{2}(x)$ is a polynomial
 of the order $P{\leq}Q-1.$ By repeating the procedure we have  $M_{2}(x)=0.$ Due to (5.8) the relation $\Phi_{2}(x)=1-\overline{\Phi_{1}(a-x)}$ holds,i.e.
 $(1/R_{2})\Psi_{2}(x)=(1/R_{1})=\overline{\Psi_{2}(a-x)}$.The last relation contradicts the condition 2 of Theorem 5.6. Hence the Proposition is true.
 \end{proof}
 
\begin{Ee}\label{Eesp}Let us consider the special case of functions of the form
(5.6):
\begin{equation}\Psi_{k}(x)=x^{m_{k}}(a-x)^{n_{k}},\quad k=1,2,\quad 0{\leq}x{\leq}a.\end{equation} \emph{We assume that $m_{k}$ and $n_{k}$ are non-negative integer and}
 \begin{equation} Q=Q_{1}=m_{1}+n_{1}{\geq}Q_{2}=m_{2}+n_{2}.\end{equation}\end{Ee}
\begin{Rk}\label{Rkco} If the relations
\begin{equation} n_{1}=m_{2},\quad m_{1}=n_{2}.\end{equation}
are true,then $\Psi_{x}(x)=\overline{\Psi_{k}(a-x)}$.Hence the zeroes of  the corresponding functions  $F_{1}(z)$ and $F_{2,1}(z)$ coincide.\end{Rk}
\begin{Tm}\label{TmOrd}Let relations (5.17) and (5.18) be true.If either one or both of equalities (5.19) are not valid ,then the  order $r$ of the corresponding differential operator $L(D)$ is defined by the relation
 \begin{equation}r=max\{n_{1}-m_{2}-1,m_{1}-n_{2}-1\}{\geq}0.\end{equation}\end{Tm}

\begin{proof}.
Let us represent the differential operator $L(D)$ in the form
 $L(D)=L_{1}(D)+L_{2}(D)$,where\\
 \begin{equation} L_{1}(D)=\sum_{p+k+s=N-1}(-1)^{k+1}\Psi_{2}^{(p)}(0)\Psi_{1}^{(k)}(0)D^{s},\end{equation}
\begin{equation}L_{2}(D)=\sum_{p+k+s=N-1}\Psi_{2}^{(k)}(a)\Psi_{1}^{(p)}(a)(-1)^{(p)}D^{s}.
\end{equation}
The order of $L_{1}(D)$ is defined by the relation $r_{1}=n_{1}-m_{2}-1$. This result follows from formula (5.21) and the equalities $r=N-1-p-k,\quad p=m_{2},\quad k=m_{1}.$ In the same way we have  $r_{2}=m_{1}-n_{2}-1$. In this case we use the equalities $p=n_{2},\quad k=n_{1}.$ The theorem is proved for the case when $r_{1}{\ne}r_{2}.$ Let us consider now the case when $r=r_{1}=r_{2}$. The coefficients by $D^{r}$ in $L_{1}$ and $L_{2}$ are respectively
\begin{equation}B_{1}=(-1)^{N+1}a^{n_{1}+n_{2}}m_{2}!m_{1}!,\quad B_{2}=(-1)^{n_{2}}a^{m_{1}+m_{2}}n_{2}!n_{1}!\end{equation} The relations
$n_{1}+n_{2}=m_{1}+m_{2}$ and $m_{2}!m_{1}!{\ne}n_{2}!n_{1}!$ imply that $B_{1}+B_{2}{\ne}0.$The theorem is proved.
\end{proof}

According to Theorems 5.6 and Proposition 5.7 the following statement is true.
\begin{Cy}\label{CyN}Let the conditions of Theorem 5.10 be fulfilled.Then all the zeroes of the corresponding functions  $F_{1}(z)$ and $F_{2,1}(z)$
haven't common zeroes. If $n_{1}{\ne}m_{1}$ ,then the corresponding function $F_{1}(z)$ hasn't real zeroes and hasn't conjugate pairs of zeroes. \end{Cy}
 \begin{Ee}\label{Een}Let us consider the case when
 \begin{equation}m_{1}=m_{2}=0,\quad n_{1}{\ne}n_{2},\quad a=1.\end{equation}\end{Ee}
 In this case we have (see [2]):
 \begin{equation}F_{k}(z)=-(-i)^{n_{k}+1}\frac{d^{n_{k}}}{dx^{n_{k}}}[x^{-1}(1-\mathrm{cos}x-i\mathrm{sin}x)],
 \quad k=1,2.\end{equation}
\begin{Cy}\label{CyZ}Let the conditions (5.24) be fulfilled.Then the corresponding functions $F_{n_{1}}(z)$ and  $F_{n_{2}}(z)$ haven't common zeroes. The corresponding function $F_{1}(z)$ hasn't real zeroes and hasn't conjugate pairs of zeroes.\end{Cy}
\begin{Rk}Using relation (1.11) we can reformulate Corollaries 5.11 and 5.13 in the terms of hypergeometric function $\Phi(b,c,z)$.\end{Rk}
\begin{Ee}\label{EeJ} Let us consider the case when
 \begin{equation}m_{1}=n_{1},\quad m_{2}=n_{2},\quad n_{1}{\ne}n_{2},\quad a=2.\end{equation}\end{Ee}
 In this case we have (see [1]):
 \begin{equation}F_{k}(z)=\sqrt{\pi}\Gamma(z/2)^{-(n_{k}+1/2)}J_{(n_{k}+1/2)}(z),\end{equation}
 where $\Gamma(z)$ is Euler Gamma function, $J_{\nu}(z)$ is Bessel function. It follows from (5.27) that the zeroes   of $F_{k}(z)$  and $J_{(n_{k}+1/2)}(z)$ other than the origin coincide.
 In case (5.26) we have $F_{k}(z)=\overline{F_{k}(\overline{z})}$. Using this fact we deduce the well-known assertion (see [1],[14]):
\begin{Cy}The functions $J_{(n_{1}+1/2)}(z)$ and $J_{(n_{2}+1/2)}(z)$ haven't common zeroes other than the origin.\end{Cy}
Now  we consider the functions of the class (1.3), where $\Psi(t)$ is a
polynomial, but we don't assume that the coefficients of  $\Psi(t)$ are algebraic numbers.
\begin{Tm}\label{Tmna}Let the following conditions be fulfilled.\\
1.The functions  $\Psi_{k}(x)$ have the form (5.6),where
\begin{equation}Q=Q_{1}{\geq}Q_{2}.\end{equation}
2. \begin{equation}\Psi_{1}(x){\ne}\overline{\Psi_{2}(a-x)}\end{equation}
3.\begin{equation}\int_{0}^{a}\Psi_{k}(x)dx{\ne}0,\quad k=1,2.\end{equation} Without loss of generality we shall suppose that
\begin{equation}\int_{0}^{a}\Psi_{k}(x)dx=1,\quad k=1,2.\end{equation}
4.\begin{equation}V(u){\not\equiv}0.\end{equation}
Then  the corresponding functions $F_{1}(z)$  and $F_{2,1}(z)$ haven't common zeroes.If $\Psi_{1}(x){\ne}\overline{\Psi_{1}(a-x)}$, then the corresponding function $F_{1}(z)$ hasn't real zeroes and hasn't conjugate pairs of zeroes. \end{Tm}
\begin{proof}.
As in the proof of Theorem 5.6 we deduce equality (5.13) which contradicts condition 4 of the theorem. It proves the theorem.
\end{proof}


\begin{thebibliography}{SaL}

\bibitem{BE}
H.Bateman and A.Erdeyi,
{\it Higher Transcendental Functions},
New York, 1963.

\bibitem{BP}
Yu.A. Brychkov and A.P. Prudnikov,
{\it Integral Transforms of Generalized Functions,} Gordon
and Breach Science Publ. 1989.

\bibitem{GHKL}
I.C. Gohberg, I. Haimovici, M.A. Kaashoek, and L. Lerer,
{\it  The Bezout integral operator: Main property and underlying abstract scheme}, in: Oper. Theory:
Adv. Appl. {\bf 161}, Birkh\"auser Verlag, Basel, 2005,  225--270.

\bibitem{GH}
I.C. Gohberg and G.Heinig, 
{\it The Continual Analogue of the Resultant Operator}, Acta  Math. Sci. Hungar
{\bf 28:(3-4)} (1976), 189-209.

\bibitem{OS}
V. Olshevsky and L. Sakhnovich,
{\it An Operator Identities Approach to Bezoutiants. A General Scheme and Examples,}
Proc. of the MTNS'04 Conference, 2004.

\bibitem{Po} M.B. Porter,
{\it On the roots of the hypergeometric and Bessel's functions},
American J. {\bf 20} (1898), 193-214.

\bibitem{SaL1}
L.A. Sakhnovich, Operatorial Bezoutiant in the theory of separation of roots of entire functions, Functional Anal. Appl. {\bf 10} (1976), 45-51 (Russian).

\bibitem{SaL2}
L.A. Sakhnovich,
 {\it Integral Equations with Difference  Kernels on Finite Intervals}, Oper. Theory:
Adv. Appl. {\bf 84}, Birkh{\"a}user Verlag,  Basel-Boston, 1996.

\bibitem{SaL3}
L.A. Sakhnovich, {\it On  the  factorization  of  the transfer matrix function},   { Sov. Math. Dokl.} {
\bf 17} (1976), 203--207.

\bibitem{SaL4}
L.A. Sakhnovich, {\it Factorisation  problems  and  operator
identities}, { Uspekhi  Mat. Nauk}  { \bf 41:1 } (1986), 3-55;
English transl.  in {  Russian Math. Surveys} { \bf 41} (1986), 1-64.

\bibitem{SaL5}
{  L.A. Sakhnovich},  {\it Spectral theory of canonical
differential
systems, method of operator identities}, Oper. Theory:
Adv. Appl.
{\bf 107}, Birkh{\"a}user Verlag,  Basel-Boston, 1999.

\bibitem{Si}
K.L. Siegel, {\it Uber einige Anwendungen diophantischer Approximationen}, Abh. preus. Acad. Wiss. {\bf 1} (1929), 1-70.

\bibitem{Ti}
E.C.Titchmarsh,
{\it Introduction to the Theory of Fourier Integrals},
Oxford,  1937.

\bibitem{Wa}
G.N. Watson,  {\it A Treatise of Bessel Functions,} Cambridge
University, 1995.
\end{thebibliography}
\end{document}